\documentclass[12pt,twoside]{article}
\usepackage[english]{babel}
\usepackage[latin1]{inputenc}
\usepackage{amsmath}
\usepackage{amssymb,amsfonts}
\usepackage{graphicx}                   

\newcommand{\bH}{\mathbf{H}}

\newcommand{\bV}{\mathbf{V}}

\newcommand{\bq}{\mathbf{q}}
\newcommand{\bp}{\mathbf{p}}

\newcommand{\bz}{\mathbf{z}}
\newcommand{\bg}{\mathbf{g}}

\newcommand{\cS}{\mathcal{S}}

\newcommand{\Smp}{\tilde{\cS}}

\newcommand{\HYP}{\mathbb{H}^3}

\newcommand{\dgr}{^\circ}

\begin{document}
\pagestyle{myheadings}
\markboth{\centerline{E.~Moln\'ar~-~I.~Prok~-~J.~Szirmai}}
{Surgeries of the Gieseking hyperbolic ideal simplex manifold $\dots$ }
\title
{Surgeries of the Gieseking hyperbolic ideal simplex manifold  
}

\author{\normalsize{E.~Moln\'ar~-~I.~Prok~-~J.~Szirmai} \\
\normalsize Budapest University of Technology and \\
\normalsize Economics Institute of Mathematics, \\
\normalsize Department of Geometry \\
\normalsize H-1521 Budapest (Hungary) \\
\date{\normalsize{\today}}}

\maketitle


\begin{abstract}
In our Novi Sad conference paper (1999) we described  
Dehn type surgeries of the famous Gieseking (1912) hyperbolic ideal 
simplex manifold $\cS$, 
leading to compact fundamental domain $\cS(k)$, $k = 2, 3, \dots$ 
with singularity geodesics 
of rotation order $k$, but as later turned out with cone angle $2(k-1)/k$. We computed also the volume of $S(k)$, tending to zero 
if $k$ goes to infinity. That time we naively thought 
that we obtained orbifolds with the above surprising property.

As the reviewer of Math. Rev., Kevin P. Scannell (MR1770996 (2001g:57030)) rightly remarked, 
``this is in conflict with the well-known theorem of D.~A. Kazhdan and 
G.~A. Margulis (1968)
and with the work of  Thurston, describing the geometric convergence 
of orbifolds under large Dehn fillings".

In this paper we refresh our previous publication. Correctly, we 
obtained cone manifolds (for $k > 2$), as A.~D.~Mednykh and V.~S.~Petrov (2006) 
kindly pointed out. We complete our discussion and derive the above cone manifold 
series (Gies.1 and Gies.2) in two geometrically equivalent form, 
by the half turn symmetry of any ideal simplex. Moreover we obtain a 
second orbifold series (Gies.3 and 4), tending 
to the regular ideal simplex as the original Gieseking manifold.  
\end{abstract}
%
\newtheorem{theorem}{Theorem}[section]
\newtheorem{corollary}[theorem]{Corollary}
\newtheorem{conjecture}{Conjecture}[section]
\newtheorem{lemma}[theorem]{Lemma}
\newtheorem{exmple}[theorem]{Example}
\newtheorem{defn}[theorem]{Definition}
\newtheorem{rmrk}[theorem]{Remark}
\newenvironment{definition}{\begin{defn}\normalfont}{\end{defn}}
\newenvironment{remark}{\begin{rmrk}\normalfont}{\end{rmrk}}
\newenvironment{example}{\begin{exmple}\normalfont}{\end{exmple}}
\newenvironment{acknowledgement}{Acknowledgement}


\section{Introduction}
The famous non-orientable Gieseking manifold (1912) is the regular simplex with 
face angles 
$\pi/3$ in the Bolyai-Lobachevskian hyperbolic space $\bH^3$ with ideal 
vertices forming a so-called cusp at the absolute, equipped by face pairing 
isometries $\bz_1$, 
$\bz_2$ as horospherical glide reflections (Fig.~1). These induce one 
equivalence class of the $6$ edges and so a ball-like neighbourhood for any 
point of an edge, and so for any point of the identified simplex $\Smp$. If 
we vary the face angles $\alpha_1$, $\alpha_2$, $\alpha_3$ at the opposite 
edges of $\Smp$, then it will be no more a manifold, but for special angles 
there exists a natural parameter \ $1<k \in \mathbb{N}$, such that $\Smp(k)$ seems to 
represent a compact hyperbolic orbifold with a closed 
singularity geodesics of rotation order $k$. We shall use the Poincaré half space 
model of $\bH^3$, with the complex projective line $\mathbb{C}_{\infty}$ 
for the absolute and with
$(w,\zeta)$, $|z_1-w||w-z_2|=\zeta \cdot \zeta$ for 
an interior point of $\HYP$, ower $w \in \mathbb{C}_{\infty}$ and third coordinate 
$\zeta$ on the half circle $z_1 z_2$ (see Fig.~1.b and Fig.~5-6). 
Moreover, We compute the volume $\bV(k)$ of $\cS(k)$ as well by 
means of the Lobachevski function $\Lambda(x)$ in \thetag{2.12}. It maybe 
surprising that $\bV(k)\to0$ if $k\to\infty$ in the conflict mentioned 
in the abstract. Our result 
implies similar consequence for the double orientable cover of the Gieseking manifold, 
i.e. for the figure-eight-knot manifold examined also by Thurston \cite{T78}. 
This paper is refreshing \cite{MPS99} as byproduct 
of [6-12]. A.~D.~Mednykh and 
V.~S.~Petrov \cite{MP06} cited our \cite{MPS99, MPS99b} and noticed our deeper mistake, 
as it will be improved here. See also \cite{AM19} with N.~V.~Abrosimov, showing 
some new phenomena and interpretations as well.

After some preliminaries in Section 2 we derive our general surgery equation 
(2.9) and our previous series in \cite{MPS99} by equation 
(2.11): an orbifold for $k = 2$ (uniformly for our series Fig.~3.b) 
and cone manifolds for $k > 2$ (Fig.~3.a and figures in the corresponding sections). 
Our exact figures and Table 2.1 show the first computer results 
(agreed with A. Mednykh and J. Weeks by different implementations). 
(2.9) yields also our further analogous results in Sections 3-5 with 
Theorems 5.1-2, Remark 5.3.

The authors thank Alexander D. Mednykh and Jeffrey R. Weeks for fruitful 
discussion and friendly help, just in the COVID-19 pandemic time.
\section{Gieseking Manifold and its Surgeries}
We start with the ideal simplex of $\bH^3$ in the half-space model, where its 
ideal vertices at infinity are represented by
\begin{equation}
0,1,z,\infty \in \mathbb{C} \cup \{\infty\}=:\mathbb{C}_\infty \tag{2.1}
\end{equation}
of the complex projective line. This will be an identified ideal simplex 
$\Smp$ with face pairing mappings $\bz_1$ and $\bz_2$ as ``horospherical 
glide reflections''
\begin{equation}
\begin{gathered}
\bz_1: \infty z \mathbf{1}=:[\bz_1^{-1}]\rightarrow \infty \mathbf{1} 
0=:[\bz_1],~{\text{i.e.}} \\
\bz_1:u \mapsto \frac{\overline{u}-1}{\overline{z}-1}~\text{or}~(u,1)\mapsto (\overline{u,1})
\begin{pmatrix}
1 & 0 \\ -1 & \overline{z}-1 \\
\end{pmatrix}, \\
\bz_2: 0 z \infty =:[\bz_2^{-1}] \rightarrow 0 \mathbf{1} z =:[\bz_2],~{\text{i.e.}} \\
\bz_2:u \mapsto \frac{\overline{u} z}{\overline{u}-\overline{z}(1-z)}~
\text{or}~(u,1)\mapsto (\overline{u,1})
\begin{pmatrix}
z & 1 \\ 0 & -\overline{z}(1-z) \\
\end{pmatrix}.
\end{gathered} \tag{2.2}
\end{equation}
As usual (e.g. in \cite{V88,V}), we extend the actions of the 
transformations into the upper half space by half-circles 
and half spheres orthogonally to the boundary plane, represented by $\mathbb{C}_\infty$.
Circles and spheres through the infinity $\infty$ will 
be orthogonal half lines and half planes, respectively. 
Thus we can describe the lines and the planes of the model space of $\bH^3$, moreover its 
congruence transformations.

Going round e.g. the edge $\overrightarrow{\infty z}$ from the starting 
identity simplex, we meet first the face $z_1^{-1}$, then follows, on the other 
side, the image face $[z_1]^{\bz_1^{-1}}$ \ $(=z_1^{-1})$ at the edge 
$(\infty 1)^{\bz_1^{-1}}$ in the $\bz_1^{-1}$-image simplex. Then the image face 
$[z_1]^{\bz_1^{-1}}$ and, on the other side, the face $[z_1]^{\bz_1^{-1} \bz_1^{-1}}$ \ 
$(=[z_1^{-1}]^{\bz_1^{-1}})$ come at edge $(\infty 0)^{\bz_1^{-1}}$ in the 
$\bz_1^{-1} \bz_1^{-1}$-image simplex. Then we meet the image face 
$[z_2^{-1}]^{\bz_1^{-1} \bz_1^{-1}}$ 
and, on the other side, the image simplex $\bz_2^{-1} \bz_1^{-1} \bz_1^{-1}$ by the 
conjugate 
mapping $\bz_1 \bz_1 \bz_2^{-1} \bz_1^{-1} \bz_1^{-1}$ of $\bz_2^{-1}$. Thus \cite{M92}, we obtain 
the cycle transformation $\bz_2 \bz_1 \bz_2^{-1} \bz_2^{-1} \bz_1^{-1} \bz_1^{-1}$ 
and we prescribe 
the trivial rotation order $\nu=1$ for the unique edge class containing $6$ 
edges. Finally we get the cycle relation 
\begin{equation}
\bz_1 \bz_1 \bz_2 \bz_2 \bz_1^{-1} \bz_2^{-1} =\mathbf{1} \tag{2.3}
\end{equation}
in equivalent form, in conformity with the fact that the dihedral angles of 
a regular ideal simplex are $\pi/3$, \ $6\cdot(\pi/3)=2\pi$ will guarantee 
ball-like neighbourhood of any point at simplex edges. However, the relation 
\thetag{2.3} with \thetag{2.2} - by careful computations - leads to equation 
\begin{equation}
|z-1|^2=|z| \tag{2.4}
\end{equation}
with more general ideal simplex, not necessarily the regular one.

Now, we turn to the ideal vertex class forming a cusp (Fig.~2--3). This 
will be represented by gluing; corresponding to images of the $4$ vertex domains 
to that of $\infty$. The side face pairing of $\Smp$ induces the pairing of 
the sides of a $2$-dimensional polygon, denoted by $\tilde{s}$ in 
Fig.~2--3, say, on a horosphere centred in $\infty$. This is represented 
in our half-space model by a Euclidean plane parallel to the absolute, and 
it can also be described on the absolute by $\mathbb{C}_\infty$.
\begin{figure}[ht]
\centering
\includegraphics[width=13cm]{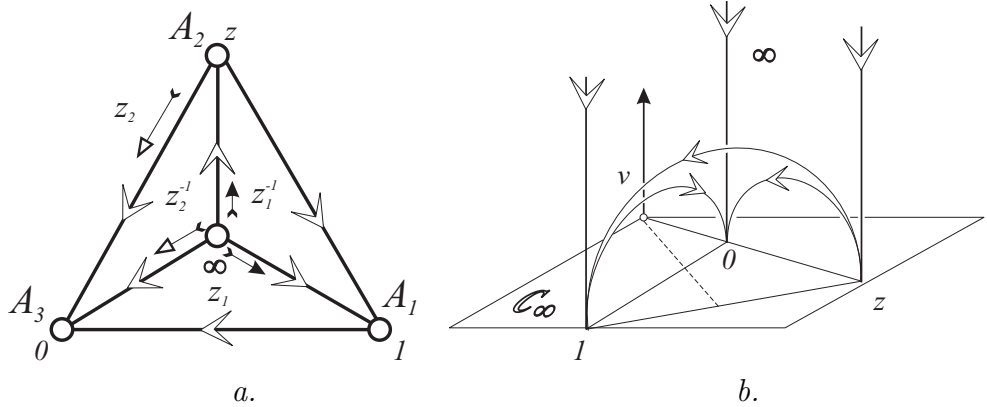}

{\it a.} \hspace{6cm} {\it b.}
\caption{The Gieseking simplex: {\it a.} with its Schlegel diagram; {\it b.} in half space model }
\label{}
\end{figure}
\begin{figure}[ht]
\centering
\includegraphics[width=11cm]{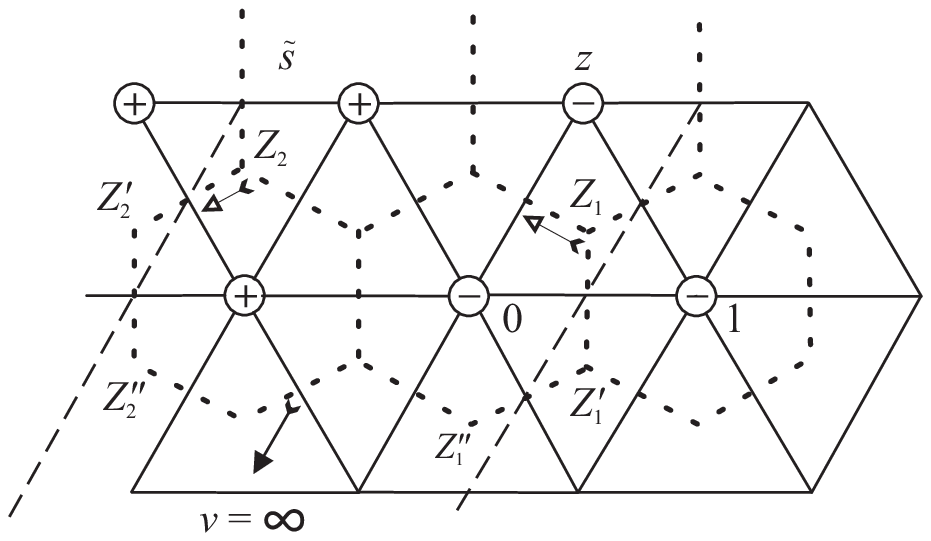}
\caption{The Gieseking regular ideal simplex tiling in the ``touching plane" at 
$\infty$, as $\mathbb{C}_\infty$}
\label{}
\end{figure}
Topologically, the polygon $\tilde{s}$ is a Klein-bottle with fundamental 
group equivariant to the Euclidean crystallographic plane group $4.~\bp \bg$. This group, 
as the stabilizer $\mathcal{G}_\infty$ of $\infty$, is determined by the starting 
group $\mathcal{G}(\bz_1,\bz_2)$ in formulae \thetag{2.2}. Fig.~3 exactly (for $k=3$
and $k=2$, respectively) 
shows the more general situation that 
$\mathcal{G}_\infty$ is generated by pairing of $s$:
\begin{equation}
\begin{gathered}
\bz_1: [z_1^{-1}]\rightarrow [z_1]~\text{``glide reflection" as before; then} \\
\bp: [z_1^{-1}]^* := [z_1^{-1}]^{\bz_2^{-1}} \rightarrow [z_1]^{*}:= 
[z_1]^{\bz_2^{-1} \bz_2^{-1}} \\
\text{i.e.}~\bp= \bz_2 \bz_1 \bz_2^{-1}  \bz_2^{-1}~=
\bz_1 \bz_1: u \rightarrow (u-z)/|z| \\
\text{a ``translation" as a central similarity in}~\mathbb{C}_\infty, \\
\bz_2^*: [z_2^{-1}]^*:= [z_2^{-1}]^{z_1^{-1} \bz_2^{-1} \bz_2^{-1}} 
\rightarrow [z_2]^{*}:= [z_2]^{\bz_1^{-1} \bz_2^{-1} \bz_2^{-1}},\\ 
{\text{i.e.}}~\bz_2^*= \bz_2 \bz_2 \bz_1 \bz_2 \bz_1^{-1} \bz_2^{-1} \bz_2^{-1}~ 
\text{a ``glide reflection"}
\end{gathered} \tag{2.5}
\end{equation}
again, it is conjugated to $\bz_2$. We see that $\bp$ is a ``translation'', 
it is $\bz_2$-conjugated to $\bz_1 \bz_2^{-1}$. 
This group $\mathcal{G}_\infty$ is 
$4.~\bp\bg$ itself (on the Euclidean plane represented by $\mathbb{C}_\infty)$ if 
$z=\frac12+i\frac{\sqrt3}2$. Then Fig.~2 shows the exact situation. We 
have obtained the Gieseking manifold with one cusp.
Other $z$, as a complex parameter, makes the stabilizer 
$\mathcal{G}_\infty$ to a 
conformal group with fixed points
\begin{equation}
\infty ~\text{and} ~ v=\frac{z}{1-|z|}. \notag
\end{equation}
\begin{figure}
\centering
\includegraphics[width=12cm]{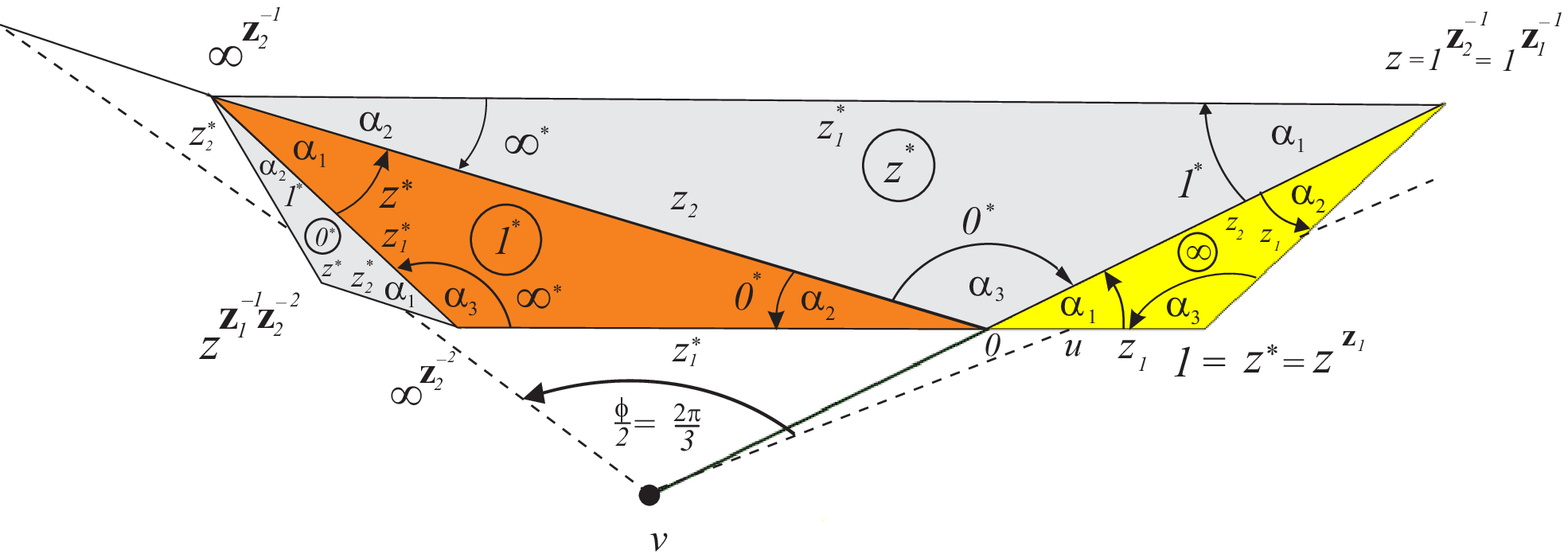} \\
\includegraphics[width=10cm]{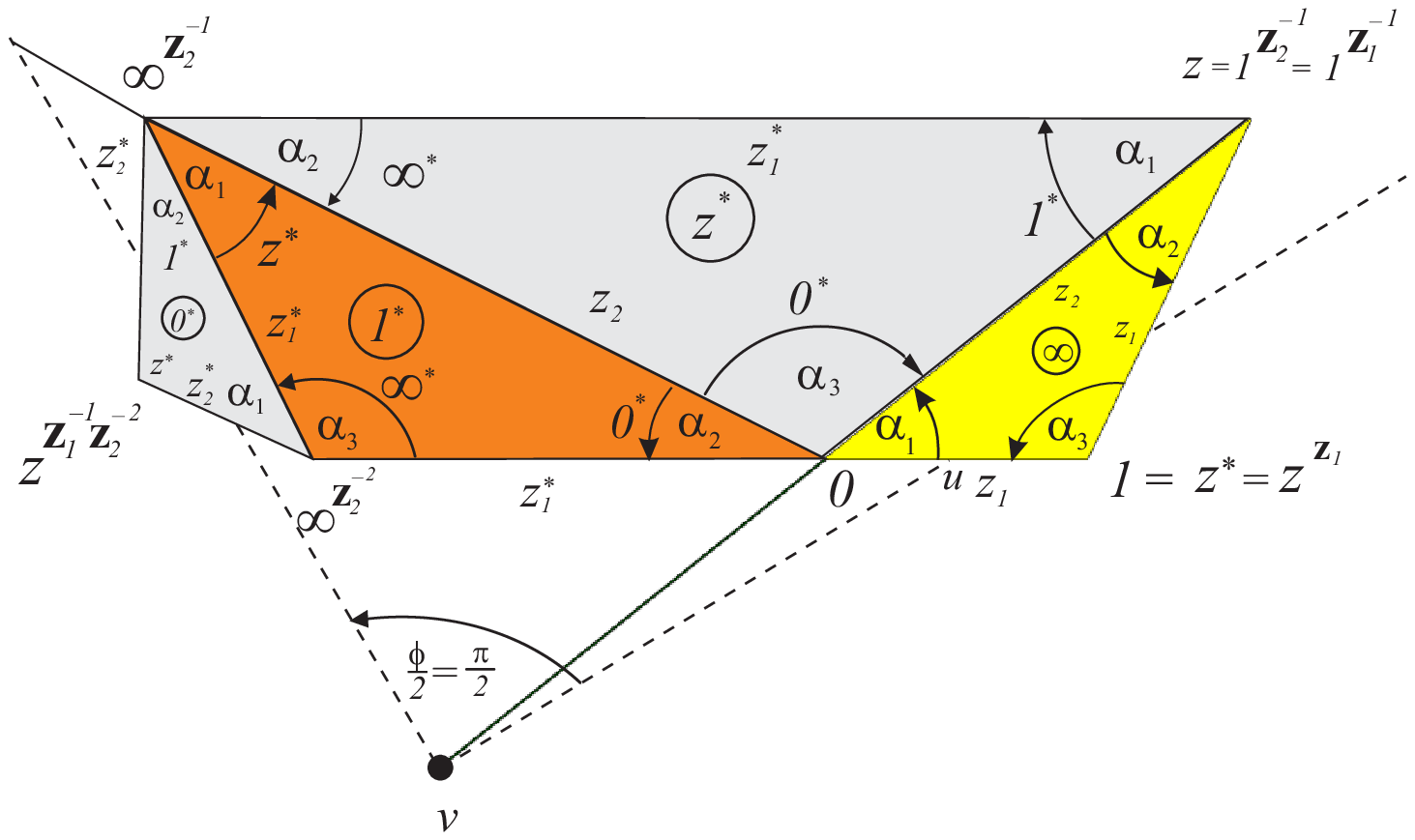}
\caption{{\it a}. The topological Klein 
bottle group $4.~\bp \bg$ in $\infty$ of $\mathbb{C}_\infty$ glued 
by fundamental domain $\mathcal{S}$; $\bz_1 \bz_2^*$ is a rotation through $-2\pi/k
=2\pi(k-1)/k$ ($\text{mod}~2\pi$) for $k=3$; {\it b}. The orbifold for $k=2$}
\label{}
\end{figure}
This line $v\infty$ will not be covered by the $\mathcal{G}_\infty$-images of the 
simplex $\Smp$ in $\bH^3$. In the model half-space the translations of 
$\mathcal{G}_\infty$ in \thetag{2.5} by \thetag{2.2} will be 
similarities of $\mathbb{C}_\infty$ with fixed 
points $v$, $\infty$. E.g. $\bz_1$ in \thetag{2.2} and $\bz_2$ in 
\thetag{2.5} are similarity-reflections indicated in Fig.~2--5. 
The simple ratio on $01$ is $u = 1/(1+|z-1|)$. For $\bz_2^*$ 
in (2.5) we can write by \thetag{2.2}
\begin{equation}
\begin{gathered}
\bz_2^*: (u,1) \rightarrow (\overline{u,1})
\begin{pmatrix}
1 & 1 \\ 0 & |1-z|^2 \\
\end{pmatrix} 
\begin{pmatrix}
1 & 0 \\ -1 & \overline{z}-1 \\
\end{pmatrix} \cdot \\
\cdot \begin{pmatrix}
\overline{z} & 1 \\ 0 & -z (1-\overline{z})
\end{pmatrix} 
\begin{pmatrix}
z-1 & 0 \\ 1 & 1 \\
\end{pmatrix}
\begin{pmatrix}
|1-z|^2 & -1 \\ 0 & 1 \\
\end{pmatrix}= \\
(\overline{u,1})
\stackrel{(2.4)}{=} (\overline{u,1})
\begin{pmatrix}
z(1-\overline{z})^2 & 0 \\ |z|^2\{|z|(\overline{z}-1)-(|z|+1)\} & -\overline{z}(1-z) \\
\end{pmatrix}.
\end{gathered} \tag{2.6}
\end{equation}
\begin{figure}
\centering
\includegraphics[width=11cm]{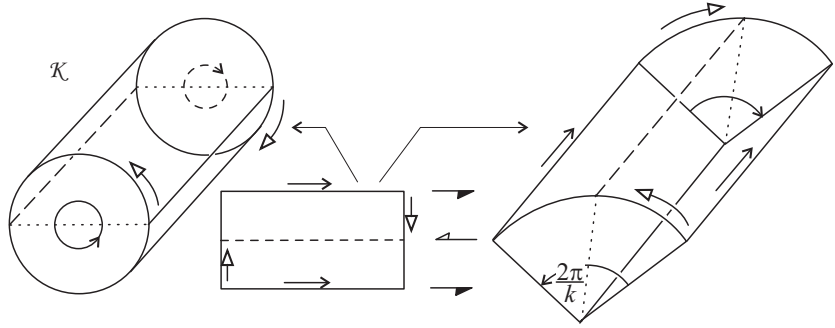}
\caption{Klein bottle solids for the later manifold or cone manifold (orbifold if the cone angle is $2\pi/k$)}
\label{}
\end{figure}
Now, we turn to the critical, so-called surgery transform $\bz_1\bz_2^*$. 
By the tricky use of (2.4),
as $$|z|=|z-1|^2=(z-1)(\overline{z}-1)=|z|^2-z-\overline{z}+1,$$ we obtain
\begin{equation}
\bz_1 \bz_2^* :(u,1) \to (u,1)
 \begin{pmatrix}
 z(1-\overline{z})^2 & 0 \\ \{|z|+1\}[|z|^2-z^2] & \overline{z}(1-z)^2 \\
 \end{pmatrix},
\tag{2.7}
\end{equation}
fixing $\infty$ and $v$, of course. We see by \thetag{2.4} that 
$\bz_1 \bz_2^*$ describes a rotation of the model half-space about the 
line $\infty v$ with angle
\begin{equation}
\phi:=\arg\left[\frac{z(1-\overline{z})^2}{\overline{z}(1-z)^2}\right]=
  2\arg z-4\arg(1-z)~(\text{mod}~2\pi),
\tag{2.8}
\end{equation}
i.e. $\phi=2 \cdot \alpha_1+4 \cdot \alpha_3$ $(\text{mod}~2\pi)$.

If we require the stabilizer $\mathcal{G}_\infty$ to act discontinuously on the 
model half-space, then this angle necessarily will be 
$\pm 2\pi/k~(\text{mod}~2\pi)$, i.e.
\begin{equation}
\begin{gathered}
\frac{z}{(1-z)^2}=\pm e^{\pm i\pi/k}\quad\text{with}\quad1<|z-1|<|z|\\
 \text{or}~ |z|<|z-1|< 1 ~\text{and}~ \mathrm{Im} z>0,\quad k=2,3,\dots
\end{gathered} \tag{2.9}
\end{equation}
\begin{figure}
\centering
\includegraphics[width=9cm]{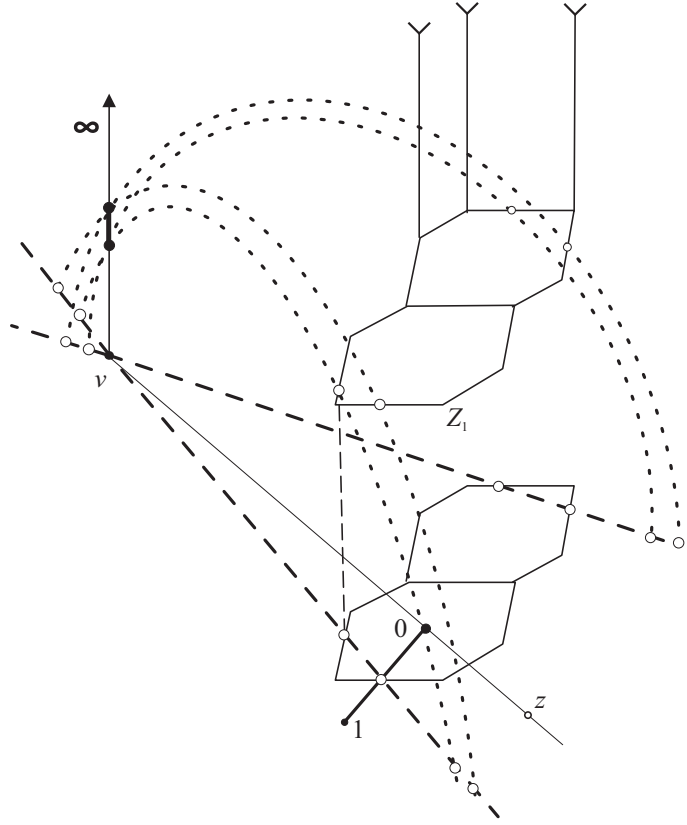}
\caption{Sketchy construction of compact fundamental domain}
\label{}
\end{figure}
can be assumed. Then $k$ is the periodicity of the rotation 
$\bz_1 \bz_2^*$ and we get 4 root series for the fundamental 
simplices, called Gies.1--Gies.4. The first one is chosen
\begin{equation}
z=1+\frac12e^{i\pi/k}\Bigl(1+\sqrt{1+4e^{-i\pi/k}}\Bigr),\qquad
k=2,3,\dots .
\tag{2.10}
\end{equation}
All data can be computed from \thetag{2.10}, especially the face angles of 
$\Smp$, equal at the opposite edges (Fig.~1.a -- 3.b)
\begin{equation}
\begin{gathered}
\left\{\begin{matrix}\overline{\infty 0}\\\overline{z1}\end{matrix}\right\}\,:\,\alpha_1=\arg z;\qquad
\left\{\begin{matrix}\overline{\infty z}\\\overline{10}\end{matrix}\right\}\,:\,\alpha_2=\arg\frac{z-1}{z};\\
\left\{\begin{matrix}\overline{\infty 1}\\\overline{z0}\end{matrix}\right\}\,:
\,\alpha_3=\arg\frac1{1-z}.
\end{gathered}
\tag{2.11}
\end{equation}
However, the computer gives more guarantees. In Table~2.1 we have computed by 
{\sc Maple} the volume of $\Smp$ as well for some values of $k$. We know 
\cite{V88,V} that the Lobachevski function
\begin{equation}
\Lambda(x)=-\int\limits_0^x\ln|2\sin\xi|d\xi\quad\text{with}\quad
\mathrm{Vol}{\Smp}=\Lambda(\alpha_1)+\Lambda(\alpha_2)+\Lambda(\alpha_3)
\tag{2.12}
\end{equation}
provides the volume of the ideal simplex with the above angles.
The formal monodromy group $\mathcal{G}(\bz_1,\bz_2,k)$ above has a unified 
``presentation"
\begin{equation}
\mathcal{G}(k)=\Bigl(\bz_1,\bz_2,\text{---},\mathbf{1}=
  \bz_1 \bz_1 \bz_2 \bz_2 \bz_1^{-1} \bz_2^{-1}=
  (\bz_1 \bz_2 \bz_2 \bz_1 \bz_2 \bz_1^{-1} \bz_2^{-2})^k_{k-1}\Bigr).
\tag{2.13}
\end{equation}
\begin{figure}[ht]
\centering
\includegraphics[width=10cm]{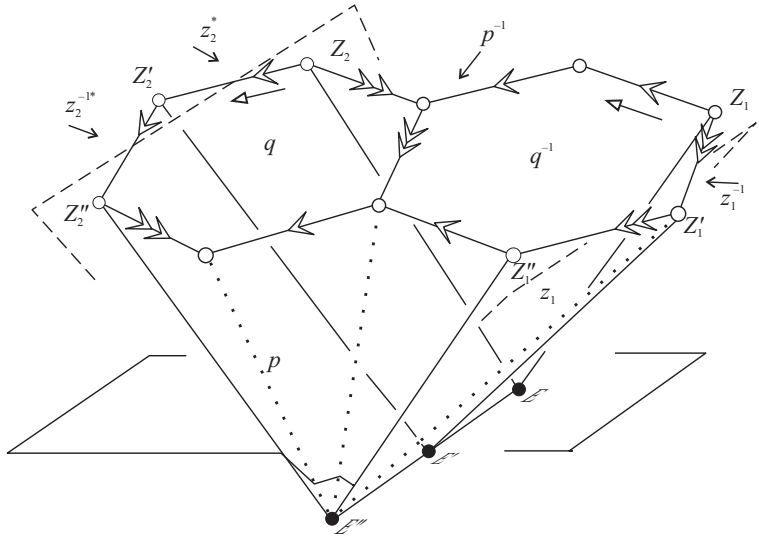}
\caption{A sketchy compact fundamental domain for $\mathcal{G}(k)$}
\label{}
\end{figure}
For $k=2$, $3$, $\dots$ we sketchily indicate by Fig.~1.b -- 6 how to construct a compact 
fundamental domain $\tilde{\mathcal{F}}_{\mathcal{G}(k)}$ (in Fig.~6) 
by deforming an ideal 
vertex domain to a compact one. We introduce an edge $(e)$ on the line 
$\infty v$ and its $\bz_1$-image with
\begin{equation}
(e)=(EE'),\qquad
e^{\bz_1}=e^{\bz_2^{* -1}}=(E'E''),\qquad
(EE') \cap (E'E'')=\emptyset. \notag
\end{equation}
Then we choose a point $Z_1$ in the simplex $\Smp=\infty01z$ and 
consider the segments
\begin{equation}
(EZ_1),\qquad
(EZ_1)^{\bz_1}=(E'Z_1'),\qquad
(EZ_1)^{\bz_1 \bz_2}=(E''Z_1''). \notag
\end{equation}
Similarly take $Z_2$, as the $\bz_1^{-1} \bz_2^{-1} \bz_2^{-1}$-image of $Z_1$ at the cusp 
gluing
\begin{equation}
(EZ_2),\qquad
(EZ_2)^{\bz_2^{* -1}}=(E'Z_2'),\qquad
(EZ_2)^{\bz_2^{* -1} \bz_2^{* -1}}=(E''Z_2''). \notag
\end{equation}
Then the corresponding curved (bent) surfaces $[\mathbf{q}^{-1}]$ and its 
$\bz_1^{-1} \bz_2^{-1} \bz_2^{-1}$-image $[\mathbf{q}]$ will be constructed, 
transversally to the edges 
of $\Smp$ (see Fig.~2,~3, 5,~6) and also \cite{M92}).
\begin{figure}
\centering
\includegraphics[width=7cm]{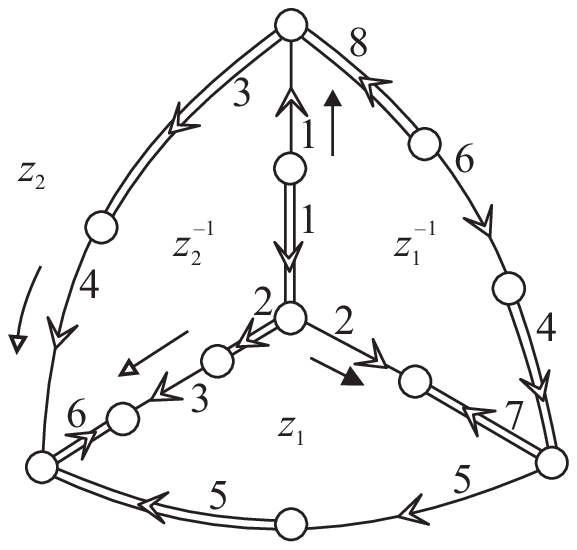}
\caption{A combinatorial fundamental domain for the Gieseking manifold 
   and cone manifold to defining relations 
   $\rightarrow~\bz_1 \bz_1 \bz_2 \bz_2 \bz_1^{-1} \bz_2^{-1}$ and 
  $\Rightarrow ~(\bz_1 \bz_2 \bz_2 \bz_1 \bz_2 \bz_1^{-1} \bz_2^{-2})^k_{k-1}$ in
  \thetag{2.13}}
\label{}
\end{figure}
Finally, in Fig.~6 we get a compact fundamental domain $\tilde{\mathcal{F}}$, with 
piecewise linear bent faces, equipped by a pairing 
$\mathcal{I}(\bz_1,\bz_2^*,\mathbf{p},\mathbf{q})$ and defining relations to 
the corresponding 
edge classes:
\begin{equation}
\bz_1\bz_1\mathbf{p}^{-1}=
\bz_2^* \bz_2^* \mathbf{p}=\bq\bz_2^{*-1} \bq^{-1} \bp \bq^{-1} \bp^{-1} =
\bz_1 \bq\bq \bp^{-1} \bq^{-1} = (\bz_1 \bz_2^*)^k=~\mathbf{1}.
\tag{2.14}
\end{equation}
In Fig.~7 we have only pictured the very economic presentation with 
its combinatorial fundamental domain:
\begin{equation}
\mathcal{G}(k)=\Bigl(\bz_1,\bz_2,~\text{---},~\mathbf{1}=
\bz_1^2 \bz_2^2 \bz_1^{-1} \bz_2^{-1} = 
(\bz_1 \bz_2^2 \bz_1 \bz_2 \bz_1^{-1} \bz_2^{-2})^k_{k-1}\Bigr).
\tag{2.15}
\end{equation}
Of course, \thetag{2.15} equivalent with \thetag{2.13} and with 
\thetag{2.14} if
\begin{equation}
\bp=\bz_1^2,\qquad
\bz_2^* =(\bz_2^2 \bz_1)  \bz_2  (\bz_1^{-1} \bz_2^{-2}),\qquad
\bq=\bz_1^{-1} \bz_2^{-1}.
\tag{2.16}
\end{equation}
\begin{rmrk}
Observe that our compactification procedure works also for $k=1$ as Fig.~4 
indicates. The cusp of our ideal simplex $\Smp$ as a Klein-bottle can be 
glued by a ``solid Klein-bottle'' $\mathcal{K}$. Then the splitting effect occurs. 
The cusp of $\Smp$ will be cut along a Klein-bottle surface to get a 
boundary. Then we glue to this boundary the boundary of $\mathcal{K}$, considered as 
$\mathbf{S}^2\!\times\!\mathbb{R}$-manifold with boundary, as follows
\begin{equation}
\mathcal{K}:=\mathbf{S}^2\!\times\!\mathbb{R}/\langle\bg\rangle,\quad
\bg:(P,r)\mapsto\Bigl(P^\mathbf{m};\;r+\frac12\Bigr),\quad
(P,r)\in \mathbf{S}^2\!\times\!\mathbb{R}.
\tag{2.17}
\end{equation}
The generator $\bg$ is a product of a reflection $\mathbf{m}$, say, in the 
equator of $\mathbf{S}^2$,  combined by a $1/2$-translation in $\mathbb{R}$. 
The boundary 
Klein-bottle can be obtained by cutting out one half-sphere, say, at 
longitudes $0$ and $\pi$, of $\mathbf{S}^2$ with complete $\mathbb{R}$-fibers.
\end{rmrk}
{\footnotesize
\smallbreak
\centerline{\vbox{
\halign{\strut\vrule~\hfil $#$ \hfil~\vrule
&\quad \hfil $#$ \hfil~\vrule
&\quad \hfil $#$ \hfil\quad\vrule
&\quad \hfil $#$ \hfil\quad\vrule
\cr
\noalign{\hrule}
\multispan4{\strut\vrule\hfill\bf Table 2.1, Cone manifold surgeries of Gieseking manifold \hfill\vrule}%
\cr
\noalign{\hrule}
   k/(k-1) &     z & Angles & Volume \cr 
\noalign{\hrule}
\noalign{\vskip2pt}
\noalign {\hrule}
   1/2&       1.624\,810\,533\,844&0.674\,888\,845\,586\approx 38.67\dgr&\cr
  i.e.~k=2 &+i\cdot1.300\,242\,590\,220&0.447\,953\,740\,604\approx 25.67\dgr&0.696\,701\,139\,104\cr
   \text{orbifold} &                           &2.018\,750\,067\,399\approx115.67\dgr&\cr
\noalign {\hrule}
   &       2.121\,964\,426\,952&0.460\,919\,465\,741\approx 26.41\dgr&\cr
  2/3&+i\cdot1.053\,755\,774\,241&0.293\,139\,042\,728\approx 16.80\dgr&0.486\,617\,604\,149\cr
   i.e.~k=3&                           &2.387\,534\,145\,121\approx136.80\dgr&\cr
\noalign {\hrule}
   &       2.327\,485\,420\,368&0.348\,223\,418\,295\approx 19.95\dgr&\cr
  3/4&+i\cdot0.844\,915\,596\,541&0.218\,587\,372\,551\approx 12.52\dgr&0.370\,676\,286\,965\cr
  i.e.~k=4 &                           &2.574\,781\,862\,743\approx147.52\dgr&\cr
\noalign {\hrule}
   &       2.558\,212\,860\,705&0.155\,851\,202\,654\approx  8.93\dgr&\cr
  8/9&+i\cdot0.401\,960\,317\,976&0.096\,607\,323\,872\approx  5.54\dgr&0.167\,339\,803\,689\cr
   i.e.~k=9&                           &2.889\,134\,127\,063\approx165.54\dgr&\cr
\noalign {\hrule}
   &       2.616\,076\,631\,698&0.017\,367\,036\,846\approx  1.00\dgr&\cr
  49/50 &+i\cdot0.073\,525\,294\,232 &0.028\,097\,779\,471\approx  1.61\dgr&0.030\,231\,732\,869\cr
  i.e.~k=50&                           &3.096\,127\,837\,270\approx177.39\dgr&\cr
\noalign {\hrule}
         &z\to(3+\sqrt5)/2\approx &\alpha_1\to0\dgr  &\cr
k \to\infty&\approx2.618\,0339      &\alpha_2\to0\dgr  &Vol \to 0\cr
         &v\to-(1+\sqrt5)/2\approx&\alpha_3\to180\dgr&\cr
         &\approx-1.618\,0339&&\cr
\noalign {\hrule} \noalign {\hrule}}}}}

\section{The second variant of our cone manifold series}
The requirements in (2.9) provide the second root series
\begin{equation}
z=1+\frac12e^{-i\pi/k}\Bigl(1-\sqrt{1+4e^{i\pi/k}}\Bigr),\qquad
k=2,3,\dots
\tag{3.1}
\end{equation}
and the cone manifold series Gies.2  for $k > 2$. 
Our Fig.~8 and Table 3.1 show these, surprisingly a little bit.

This will be geometrically equivalent to Gies.1 
by the half-turn symmetry of ideal simplex: 
$0 \leftrightarrow \infty$, $1 \leftrightarrow z$.
\begin{figure}
\centering
\includegraphics[width=10cm]{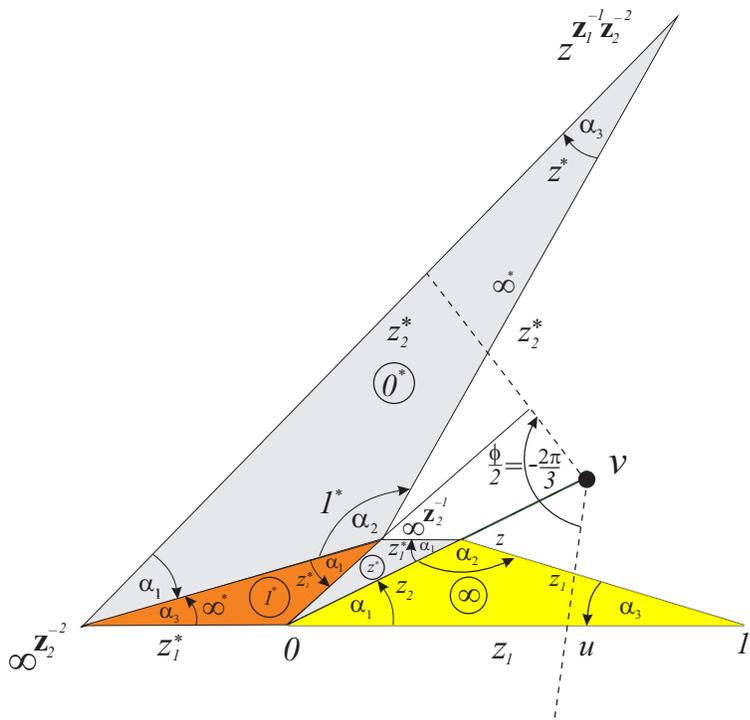}
\caption{Gies.2 series by cone manifold for $k = 3$}
\label{}
\end{figure}
{\footnotesize
\smallbreak
\centerline{\vbox{
\halign{\strut\vrule~\hfil $#$ \hfil~\vrule
&\quad \hfil $#$ \hfil~\vrule
&\quad \hfil $#$ \hfil\quad\vrule
&\quad \hfil $#$ \hfil\quad\vrule
\cr
\noalign{\hrule}
\multispan4{\strut\vrule\hfill\bf Table 3.1, Second variant of Gieseking cone manifolds  \hfill\vrule}%
\cr
\noalign{\hrule}
   (k-1)/k &     z & Angles & Volume \cr 
\noalign{\hrule}
\noalign{\vskip2pt}
\noalign {\hrule}
   &       0.375\,189\,466\,155 &0.674\,888\,845\,586\approx 38.67\dgr&\cr
  1/2 &+i\cdot 0.300\,242\,590\,218 &2.018\,750\,067\,399\approx115.67\dgr&0.696\,701\,139\,104\cr
   \text{orbifold}&                           &0.447\,953\,740\,604\approx 25.67 \dgr&\cr
\noalign {\hrule}
   &       0.378\,035\,573\,048&0.460\,919\,465\,741\approx 26.41\dgr&\cr
  2/3&+i\cdot0.187\,730\,370\,454&2.387\,534\,145\,121\approx136.80\dgr&0.486\,617\,604\,149\cr
  &                           &0.293\,139\,042\,728\approx 16.80\dgr&\cr
\noalign {\hrule}
   &       0.379\,621\,360\,824&0.348\,223\,418\,295\approx 19.95\dgr&\cr
  3/4&+i\cdot0.137\,808\,815\,354&2.574\,781\,862\,743\approx147.52\dgr&0.370\,676\,286\,965\cr
  &                           &0.218\,587\,372\,551\approx 12.52\dgr&\cr
\noalign {\hrule}
   &       0.381\,479\,760\,078&0.155\,851\,202\,654\approx  8.93\dgr&\cr
  8/9&+i\cdot0.059\,940\,174\,652&2.889\,134\,127\,063\approx165.54\dgr&0.167\,339\,803\,689\cr
  &                           &0.096\,607\,323\,872\approx  5.54\dgr&\cr
\noalign {\hrule}
   &       0.381\,950\,096\,732&0.017\,367\,036\,846\approx  1.00\dgr&\cr
  49/50 &+i\cdot0.010\,734\,774\,696 &3.096\,127\,837\,270\approx177.39\dgr&0.030\,231\,732\,869\cr
  &                           &0.028\,097\,779\,471\approx  1.61\dgr&\cr
\noalign {\hrule}
         &z\to(3-\sqrt5)/2\approx &\alpha_1\to0\dgr  &\cr
k \to\infty&\approx0.381\,966\,012  &\alpha_2\to180\dgr  &Vol \to 0\cr
         &v=\frac{z}{1-|z|} \to \frac{3-\sqrt5}{\sqrt{5}-1}\approx&\alpha_3\to0\dgr&\cr
         &\approx 0.618\,033\,989& &\cr
\noalign {\hrule} \noalign {\hrule}}}}}
\smallbreak
\section{Gies.3-4 tend to the regular ideal simplex manifold}
The requirements in (2.9) provide the third root series
\begin{equation}
z=1-\frac12e^{-i\pi/k}\Bigl(1+\sqrt{1-4e^{i\pi/k}}\Bigr),\qquad
k=2,3,\dots
\tag{4.1}
\end{equation}
and the orbifold series Gies.3. Our Fig.~9 show the case $k = 3$ and $9$.  
Table 4.1 gives computer results.
\begin{figure}
\centering
\includegraphics[width=11cm]{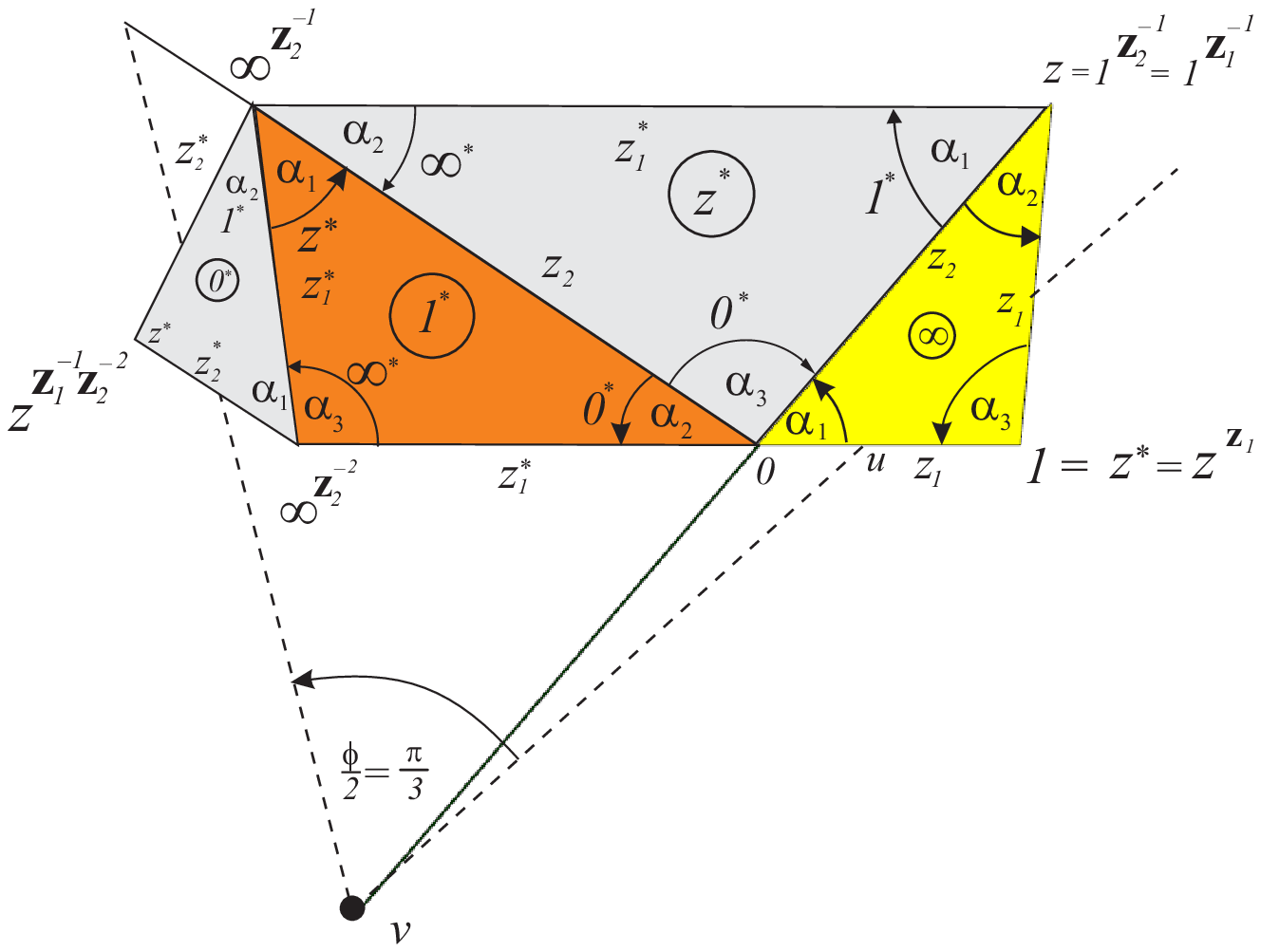}\\
\includegraphics[width=11cm]{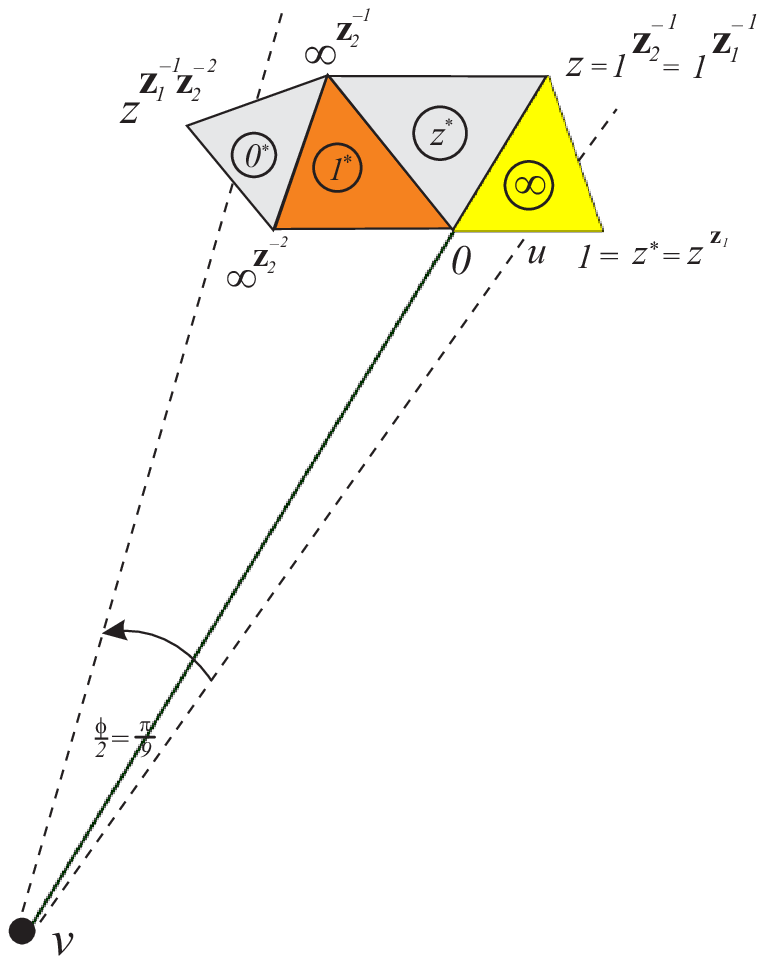}

\caption{Gies.3 series is represended by orbifolds {\it a.} $k = 3$; {\it b.} $k=9$}
\label{}
\end{figure}
The fourth root series will be
\begin{equation}
z=1-\frac12e^{i\pi/k}\Bigl(1-\sqrt{1-4e^{-i\pi/k}}\Bigr),\qquad
k=2,3,\dots
\tag{4.2}
\end{equation}
and the orbifold series Gies.4. Both last series tend to the Gieseking 
regular ideal simplex manifold, as J.~R.~Weeks predicted for us in our discussions.
{\footnotesize
\smallbreak
\centerline{\vbox{
\halign{\strut\vrule~\hfil $#$ \hfil~\vrule
&\quad \hfil $#$ \hfil~\vrule
&\quad \hfil $#$ \hfil\quad\vrule
&\quad \hfil $#$ \hfil\quad\vrule
\cr
\noalign{\hrule}
\multispan4{\strut\vrule\hfill\bf Table 4.1, Gies.3 orbifolds tending to Gieseking manifold  \hfill\vrule}%
\cr
\noalign{\hrule}
   k &     z & Angles & Volume \cr 
\noalign{\hrule}
\noalign{\vskip2pt}
\noalign {\hrule}
   &       1.624\,810\,533\,840 & 0.674\,888\,845\,586\approx 38.67\dgr&\cr
  k=2 &+i\cdot 1.300\,242\,590\,218 &0.447\,953\,740\,604\approx 25.67\dgr&0.696\,701\,139\,104\cr
   &                           &2.018\,750\,067\,399\approx115.67\dgr&\cr
\noalign {\hrule}
   &       1.121\,744\,414\,125&0.861\,384\,224\,935\approx 49.35\dgr&\cr
  k=3&+i\cdot1.306\,622\,402\,750&0.616\,505\,438\,729\approx 35.33\dgr&0.865\,129\,197\,896\cr
  &                           &1.663\,702\,989\,926\approx 95.32\dgr&\cr
\noalign {\hrule}
   &       0.906\,141\,518\,670&0.937\,270\,974\,150\approx 53.70\dgr&\cr
  k=4&+i\cdot1.233\,637\,606\,497&0.709\,461\,758\,021 \approx 40.65\dgr &0.928\,730\,994\,975\cr
  &                           &1.494\,859\,921\,419\approx 85.65\dgr&\cr
\noalign {\hrule}
   &       0.634\,217\,447\,865&1.024\,132\,633\,920\approx 58.68\dgr&\cr
  k=9&+i\cdot1.042\,222\,316\,635&0.884\,134\,127\,063\approx50.66\dgr&0.997\,471\,628\,531\cr
  &                           &1.233\,262\,935\,034\approx  70.66\dgr&\cr
\noalign {\hrule}
   &       0.519\,144\,497\,912&1.046\,438\,204\,847\approx 59.96\dgr&\cr
  k=50 &+i\cdot0.897\,609\,875\,813 &1.016\,161\,297\,835\approx58.22\dgr&1.014\,371\,909\,442\cr
  &                           &1.078\,993\,150\,908\approx  61.82\dgr&\cr
\noalign {\hrule}
  &       & \alpha_1 \to 60\dgr=\frac{\pi}{3} \approx 1.047\,197\,551\,196 &\cr
  k \to\infty& z\to \frac{1}{2} + i\cdot\frac{\sqrt{3}}{2}&\alpha_2 \to 60\dgr=\frac{\pi}{3} \approx 1.047\,197\,551\,196&1.014\,941\,606\,409\cr
  &         v \to \infty                  &\alpha_3 \to 60\dgr=\frac{\pi}{3} \approx 1.047\,197\,551\,196&\cr
\noalign {\hrule} \noalign {\hrule}}}}}
We do not give here illustration to Gies.4 series, 
equivalent to the previous one by half-turn symmetry again.

As we see in Table 4.1 our orbifold volumes 
(more important, the half of them) are large enough. T.~H.~Marshall and G.~J.~Martin
\cite{MM12} determined the exact lower bound $\approx 0.0390$, 
the next is $\approx 0.0408$ for orientable orbifolds in dimension three. 
Compare that the half of the Coxeter orthoscheme $(5, 3, 5)$  $\approx 0.0467$ and 
$(3, 5, 3)$ $\approx 0.01953$, two-times less than 
the optimal one, but this 
orthoscheme has also reflection, as the authors noticed as well. 
In higher dimensions the problem is open, in general.
\section{Summary}
Now we summarize our results.

\begin{theorem}
The surgery procedure of Gieseking manifold leads essentially to two 
different series. 
For rotation parameter $k = 2$ we get an orbifold.  
For $k>2$ the surgery yields compact 
nonorientable hyperbolic cone manifolds in the first case Gies.1-2, 
with underlying Gieseking manifold 
before, where a closed geodesic line exists with cone angle $2\pi(k-1)/k$.
This can be realized by a deformed ideal simplex $\mathcal{S}(k)$ by 
\thetag{2.1--2} with complex parameter $z(k)$ by \thetag{2.11} that 
uniquely determines all metric data in our figures and tables. The volume of 
$\mathcal{S}(k)$ tends to $0$ if $k\to\infty$. $\square$
\end{theorem}

\begin{theorem}
The second series in cases Gies.3-4 leads to orbifolds also for 
$k \ge 2$, tending to the original Gieseking manifold if $k$ goes to 
infinity. 
\end{theorem}
\begin{rmrk}
The orientable double cover of Gieseking manifold, known as 
Thurston manifold (or the complement of the figure-eight-knot) 
has a ``manifold surgery" of volume $0.9813688289 \dots$, 
which is known as second minimal one. But the above construction 
leads to cone manifold surgeries whose volumes tend to zero. 
to the first series and tend to the original manifold to the second series.

The minimal volume orientable manifold of Fomenko-Matveev-Weeks 
(with volume $0.94270736 \dots$) can also be obtained by surgery. 
The occasionally possible (?) orbifold surgery is not examined yet (?).

We found in \cite{MPS99b} the third non-orientable 
double-ideal-regular-simplex-ma\-ni\-fold
by computer. This has similar surgery phenomena 
that will be published again, because of its actuality.
We plan to discuss the general surgery situations by refreshing 
\cite{MPS99b} as well.

\end{rmrk}


\noindent
\footnotesize{Budapest University of Technology and Economics Institute of Mathematics, \\
Department of Geometry, \\
H-1521 Budapest, Hungary. \\
E-mail:emolnar@math.bme.hu,~prok@math.bme.hu,~szirmai@math.bme.hu \\


\begin{thebibliography}{999999999}

%
\bibitem{AM19} N.~V.~Abrosimov, A.~D.~Mednykh, 
Area and volume in non-Euclidean geometry, Chapter 11 in the book
Eighteen Essays in non-Euclidean geometry,
\textit{European Mathematical Society Publishing House, Monograph series, 2019,
ISBN print 978-3-03719-196-5, ISBN online 978-3-03719-696-0},
151-189.
%
\bibitem{FM88} A.~T.~Fomenko, S.~V.~Matveev,
Isoenergetic surfaces of Hamiltonian systems,
account of three-dimensional manifolds in order of their complexity
and computation of volumes of closed hyperbolic manifolds,
\textit{Uspehi~mat. nauk}, {\bf{43}} (1988), \rm 5--22, (Russian).
%
\bibitem{KM68} D.~A.~Kazhdan, G.~A.~Margulis,
A proof of Selberg's hypothesis,
\textit{Mat.Sb. (N.S.)}, {\bf{75/117}} (1968), \rm 163--168, (Russian).
%
%
\bibitem{MM12} T.~H.~Marshall, G.~J.~Martin,
Minimal co-volume hyperbolic lattices, II: Simple torsion in a Kleinian group,
\textit{Annals of Mathematics}, {\bf{176}} (2012), \rm 261--301.
%
\bibitem{MP06} A.~D.~Mednykh, V.~S.~Petrov, 
On spontaneous Surgery on Knots and Links, 
{\it Non-Euclidean Geometries, J\'{a}nos Bolyai Memorial 
Volume, Editors: A. Pr\'{e}kopa and E. Moln\'{a}r, 
Mathematics and Its Applications}, {\bf 581}, Springer (2006), 307--319.%
%
\bibitem{M92} E.~Moln{\'a}r,
Polyhedron complexes with simply transitive group actions
and their realizations,
\textit{Acta Math. Hung.}, {\bf{59(1-2)}} (1982), \rm 175--216.
%
\bibitem{MP91} E.~Moln{\'a}r, I.~Prok,
Classification of solid transitive simplex tilings in simply
  connected 3-spaces,
  Part~{I.} The combinatorial description by figures and tables,
\textit{Colloquia Math. Soc. J{\'a}nos Bolyai
  Intuitive Geometry,
  {\rm Szeged (Hungary), 1991}
North--Holland Publ. Comp.
Amsterdam--Oxford--New~York}, {{\bf63.}} (1994), \rm 311--362.
%
\bibitem{MPS97} E.~Moln{\'a}r, I.~Prok, J.~Szirmai,
Classification of solid transitive simplex tilings in simply  connected
  3-spaces, Part~{II.} Metric realizations of the maximal simplex tilings,
\textit{Periodica Math. Hung.}, {\bf{35(1-2)}} (1997), \rm 47--94.
%
\bibitem{MPS99} E.~Moln{\'a}r, I.~Prok, J.~Szirmai,
The Gieseking manifold and its surgery orbifolds,
\textit{Novi Sad, Journal of Mathematics}, {\bf{29(3)}} (1999), \rm 187--197.
%
\bibitem{MPS99b} E.~Moln{\'a}r, I.~Prok, J.~Szirmai,
Classification of hyperbolic manifolds and related orbifolds with 
charts up to two ideal simplices,
\textit{Proceedings of "Internationale Tagung über Geometrie, 
Algebra und Analysis" Balatonfüred, Hungary}, (1999), \rm 293--315.
%
\bibitem{MPS06} E.~Moln\'{a}r, I.~Prok, J.~Szirmai, 
Classification of tile-transitive 3-simplex tilings and their 
realizations in homogeneous spaces, 
{\it Non-Euclidean Geometries, J\'{a}nos Bolyai Memorial 
Volume, Editors: A. Pr\'{e}kopa and E. Moln\'{a}r, 
Mathematics and Its Applications}, {\bf 581}, Springer (2006), 321--363.
%
\bibitem{P92} I.~Prok,
Data structures and procedures for a polyhedron algorithm,
\textit{Periodica Polytechnica Ser. Mech. Eng.}, {\bf{36(3-4)}} (1992), \rm 299-316.
%
\bibitem{T78} W.~Thurston,
The geometry and topology of 3-manifolds,
\textit{Princeton University (Lecture notes)}, 1978.
%
\bibitem{V88} E.~B.~Vinberg, O.~V.~Shvartsman,
Discrete transformation groups of spaces of constant curvature,
\textit{Geometriya~2 VINITI Itogi Nauki i Techniki, 
  Sovr. Probl. Mat. Fund. Napr.}, \bf{29} (1988), \rm 147--259, (Russian).
%
\bibitem{V}
E.~B.~Vinberg (Ed.),
\textit{Geometry II. Spaces of Constant Curvature}
Spriger Verlag Berlin-Heidelberg-New York-London-Paris-Tokyo-Hong Kong-Barcelona-Budapest, 1993.
%
\bibitem{W85}
J.~R.~Weeks,
\textit{Hyperbolic structures on three-manifolds}
PhD dissertation, Princeton 1985.

\end{thebibliography}
\end{document}